\documentclass[dvips, 10pt, a4paper]{article}
\usepackage{latexsym, amsmath}
\usepackage{amssymb}

\begin{document}

\centerline{\Large{Generalized inequalities on warped product submanifolds}}

\smallskip
\centerline{\Large{in nearly trans-Sasakian manifolds}}

\bigskip

\centerline{{\it{Abdulqader Mustafa, Siraj Uddin and B.R. Wong}}}

\vspace*{.15cm}
\begin{abstract}
In this paper, we study warped product submanifolds of nearly trans-Sasakian manifolds. The non-existence of the warped product semi-slant submanifolds of the type $N_\theta\times{_{f}N_T}$ is shown, whereas some characterization and new geometric obstructions are obtained for the warped products of the type  $N_T\times{_{f}N_\theta}$. We establish two general inequalities for the squared norm of the second fundamental form. The first inequality generalizes derived inequalities for some contact metric manifolds [16, 18, 19, 24], while by a new technique, the second inequality is constructed to express the relation between extrinsic invariant (second fundamental form) and intrinsic invariant (scalar curvatures). The equality cases are also discussed.\\

\end{abstract}

\section{Introduction}

In a natural way, warped products appeared in differential geometry generalizing the class of Riemannian product manifolds to a much larger one, called warped product manifolds, which are applied in general relativity to model the standard space time, especially in the neighborhood of massive stars and black holes [20, 21]. These manifolds were introduced by Bishop and O'Neill [3]. They defined warped products as follows: Let $N_1$ and $N_2$ be two Riemannian manifolds with Riemannian metrics $g_1$ and $g_2$, respectively, and $f>0$ be a differential function on $N_1$. Consider the product manifold $N_1\times N_2$ with its projections $\pi_1:N_1\times N_2\rightarrow N_1$ and $\pi_2:N_1\times N_2\rightarrow N_2$. Then their warped product manifold $M= N_1\times _fN_2$ is the Riemannian manifold $N_1\times N_2=(N_1\times N_2, g)$ equipped with the Riemannian structure such that
$$\|X\|^2= \|{\pi_1}_\star (X)\|^2+ (f\circ\pi_1)^2 \|{\pi_2}_\star(X)\|^2,$$
for any vector field $X$ tangent to $M$, where $\star$ is the symbol for the tangent maps. A warped product manifold $M=N_1\times N_2$ is said to be {\it{trivial}} or simply {\it{Riemannian product}} if the warping function $f$ is constant. For the survey on warped products as Riemannian submanifolds we refer to [12].

\parindent=8mm
A $(2m+1)-$dimensional $C^\infty$ manifold $(\bar M, g, \phi, \xi, \eta)$ is said to have an {\it{almost contact structure}} if there exist on $\bar M$ a tensor field $\phi$ of type $(1,1)$, a vector field $\xi$, a $1-$form $\eta$ and a Riemannian metric $g$ satisfying

$$\phi^2=-I+\eta\otimes\xi,~~~~~\phi\xi=0,~~~~~\eta\circ\phi=0,~~~~~\eta(\xi) = 1,\eqno(1.1)$$
$$~~~~~\eta(X)=g(X, \xi),~~~~~g(\phi X, \phi Y)=g(X, Y)-\eta(X)\eta(Y).\eqno(1.2)$$
where $X$ and $Y$ are vector fields on $\bar M$[5]. We shall use the symbol $\Gamma(T\bar M)$ to denote the Lie algebra of vector fields on the manifold $\bar M$.

\parindent=8mm
In the classification of almost contact structures, D. Chinea and C. Gonzalez [14] divided this structure into twelve well known classes, one of the classes that appears in this classification is denoted $C_1\oplus C_5 \oplus C_6$, according to their classification an almost contact metric manifold is a nearly trans-Sasakian manifold if it belongs to this class. Another line of thought C. Gherghe introduced nearly trans-Sasakian structure of type $(\alpha , \beta )$, which generalizes trans-Sasakian structure in the same sense as nearly Sasakian generalizes Sasakian ones, in this sense an almost contact metric structure $(\phi , \xi , \eta , g )$ on $\bar M$ is called a {\it{nearly trans-Sasakian structure}} if
$$(\bar\nabla_X\phi)Y+(\bar\nabla_Y\phi)X=\alpha (2g(X, Y)\xi -\eta(Y)X-\eta(X)Y)~~~~~~~~~~~~~$$
$$~~~-\beta(\eta(Y)\phi X+\eta(X)\phi Y)\eqno(1.3)$$
for any $X, Y\in\Gamma(T\bar M)$. Moreover, a nearly trans-Sasakian of type $(\alpha , \beta)$ is nearly-Sasakian, or nearly Kenmotsu, or nearly cosymplectic according as $\beta$ = 0 or $\alpha$ = 0 or $\alpha=\beta=0$.

\parindent=8mm
J. S. Kim et.al [17] initiated the study of semi-invariant submanifolds of nearly trans-Sasakian manifolds and obtained many results on the extrinsic geometric aspects of these submanifolds, whereas the slant submanifolds were studied in the setting of nearly trans-Sasakian manifolds by Al-Solamy and V.A. Khan [1]. Recently, we have initiated the study of CR-warped product in nearly trans-Sasakian manifolds [19]. In the present paper we consider warped product of proper slant and invariant submanifolds of nearly trans-Sasakian manifolds, called warped product semi-slant submanifolds. The paper is organized as follows: Section 2 is devoted to provide the basic definitions and formulas which are useful to the next section. In section 3, general and special non-existence results are proved. In section 4, the necessary lemmas for the two inequalities and some geometric obstructions are obtained. In section 5, a general inequality which generalizes obtained inequalities in [16, 18, 19, 24] is established. In section 6, we develop a new technique by means of Gauss equation and apply to construct a general inequality for the second fundamental form in terms of the scalar curvatures of submanifolds and the warping function.

\section{Preliminaries}
\parindent=8mm Let $M$ be a $n$-dimensional Riemannian manifold isometrically immersed in any Riemannian manifold $\bar M$. Then, Gauss and Weingarten formulae are respectively given by
$$\bar \nabla_X Y=\nabla_X Y+h(X,Y), \eqno(2.1)$$
and
$$\bar\nabla_XN=-A_NX+\nabla^\perp_XN, \eqno(2.2)$$
for all $X,Y\in\Gamma(TM)$, where $\nabla$ is the induced
Riemannian connection on $M$, $N$ is a vector field normal to $\bar
M$, $h$ is the second fundamental form of $M$, $\nabla^\perp$ is the
normal connection in the normal bundle $T^\perp M$ and $A_N$ is the
shape operator of the second fundamental form. They are related as
$$g(A_NX,Y)=g(h(X,Y),N),\eqno(2.3)$$
where $g$ denotes the Riemannian metric on $\bar M$ as
well as the metric induced on $M$. For any $X\in\Gamma(TM)$, we decompose $\phi X$ as follows
$$\phi X=PX+FX\eqno(2.4)$$
where $PX$ and $FX$ are the tangential and normal components of $\phi X$, respectively.

\parindent=8mm For a submanifold $M$ of an almost contact manifold $\bar M$, if $F$ is identically zero then $M$ is $invariant$ and if $P$ is identically zero then $M$ is $anti-invariant$.

\parindent=8mm
For any orthonormal basis $\{e_1, \cdots , e_n\}$ of the tangent space $T_xM$, the $\it mean~ curvature~ vector$ $\vec H(x)$ is given by
$$\vec H(x)=\frac{1}{n} \sum_{i=1}^{n} h(e_i, e_i)$$
where $n=dim (M)$. The submanifold $M$ is {\it totally geodesic} in $\bar M$ if $h=0,$ and {\it minimal} if $H=0$. If $h(X,Y) =g(X,Y)H$ for all $X, Y \in\Gamma(TM)$, then $M$ is {\it totally umbilical}.\\

Let $(M, g)$ be a submanifold of a Riemannian manifold $\bar M$ equipped with a Riemannian metric $g$. The $\it equation~ of~ Gauss$ is given by
$$R(X,Y,Z,W)=\bar R(X,Y,Z,W)+g(h(X, W), h(Y, Z))$$
$$-g(h(X, Z), h(Y, W)),\eqno (2.5)$$
for all $X, Y, Z, W \in\Gamma(TM)$, where $\bar R$ and $R$ are the curvature tensors of $\bar M$ and $M$ respectively, and $h$ is the second fundamental form.\\

\noindent
{\bf{Definition 2.1 [11].}} {\it{An immersion $\varphi : N_1\times _fN_2\rightarrow \bar M$ is called $N_i$-{\it totally~ geodesic} if the partial second fundamental form $h_i$ vanishes identically. It is called $N_i$-{\it minimal} if the partial mean curvature vector $\vec H_i$ vanishes, for $i=1, 2$.}}

\parindent=8mm
The {\it scalar curvature} $\tau (x)$ of $M$ is defined by
$$ \tau (x)=\sum_{1\le i< j\le n} K (e_i \wedge e_j),\eqno (2.6)$$
where $K(e_i\wedge e_j)$ is the $\it sectional $ $\it curvature$ of the plane section spanned by $e_i$ and $e_j$ at $x\in M$. Let $\Pi_k$ be a $k$-plane section of $T_xM$ and $\{e_1, \cdots , e_k\}$ any orthonormal basis of $\Pi_k$. The scalar curvature $\tau (\Pi_k)$ of $\Pi_k$ is given by [11]
$$ \tau (\Pi_k)=\sum_{1\le i< j\le k} K (e_i \wedge e_j).$$
The scalar curvature of $\tau (x)$ of $M$ at $x$ is identical with the scalar curvature of the tangent space $T_xM$ of $M$ at $x$, that is, $\tau (x)= \tau (T_xM)$. Geometrically, $\tau (\Pi_k)$ is the scalar curvature of the image $\exp_x(\Pi_k)$ of $\Pi_k$ at $x$ under the exponential map at $x$. If $\Pi_2$ is a $2$-plane section, $\tau (\Pi_2)$ is simply the sectional curvature $K(\Pi_2)$ of $\Pi_2$, [11, 12, 13].

\parindent=8mm Now, let us put $$h_{ij}^r=g(h(e_i,e_j),e_r),\eqno (2.7)$$
where $i, j \in \{1,\cdots, n\}$, and $r \in \{ n+1,\cdots, 2m+1\}$. Then, in view of the equation of Gauss, we have
$$ K (e_i \wedge e_j )= \bar K (e_i \wedge e_j) + \sum_{r=n+1}^{2m+1} (h_{ii}^{r} h_{jj}^{r}- (h_{ij}^{r})^2),\eqno(2.8)$$
where  $K (e_i \wedge e_j )$ and $ \bar K (e_i \wedge e_j)$ denote to the sectional curvature of the plane section spanned by $e_i$ and $e_j$ at $x$ in the submanifold $M$ and in the ambient manifold $\bar M$ respectively. Taking the summation over the orthonormal frame of the tangent space of $M$ in the above equation, we obtain
$$2\tau (x)= 2\bar\tau (T_xM) +n^2 \|H\|^2-\|h\|^2,\eqno(2.9)$$
where $\bar\tau (T_xM)=\sum_{1\le i< j\le n}\bar K (e_i \wedge e_j)$ denotes the scalar curvature of the n-plane section $T_xM$, for each $x\in M$ in the ambient manifold $\bar M$.

\parindent=8mm There are different classes of submanifolds which we introduce briefly such as slant submanifolds, CR-submanifolds and semi-slant submanifolds. We shall always consider $\xi$ to be tangent to the submanifold $M$. For a slant submanifold $M$,  there is a non zero vector $X$ tangent to $M$ at $x$, such that $X$ is not proportional to $\xi_x$, we denote by $0\leq\theta (X) \leq\pi /2$, the angle between $\phi X$ and $T_x M$ is called the slant angle. If the slant angle $\theta (X)$ is constant for all $X\in T_x M - \langle\xi_x\rangle$ and $\ x\in M$, then M is said to be a slant submanifold [7]. Obviously, if $\theta=0$, M is invariant and if $\theta$ = $\pi /2$, M is an anti-invariant submanifold. A slant submanifold is said to be $proper~ slant$ if it is neither invariant nor anti-invariant submanifold.

\parindent=8mm
We recall the following result for a slant submanifold of an almost contact metric manifold.\\

\noindent
{\bf {Theorem 2.1 [7].}}{  Let M be a submanifold of an almost contact metric manifold $\bar M$, such that $\xi \in\Gamma(TM)$. Then M is slant if and only if there exists a constant  $\lambda\in [0, 1]$ such that
$$P^2=\lambda (- I + \eta \otimes \xi)\eqno (2.10)$$
Furthermore, if $\theta$ is slant angle, then $\lambda = \cos ^2\theta$.

\parindent=8mm
Following relations are straightforward consequence of equation (2.10)
$$g(PX,PY)=\cos^2\theta ( g(X,Y)- \eta (Y)\eta (X))\eqno(2.11)$$
$$g(FX,FY)=\sin^2\theta ( g(X,Y)-\eta (Y)\eta (X))\eqno(2.12)$$
for all $X,Y\in\Gamma(TM)$.

\parindent=8mm The idea of semi-slant submanifolds of almost Hermitian manifolds was given by N. Papaghuic [22]. In fact, semi-slant submanifolds were defined on the line of CR-submanifolds. These submanifolds are defined and investigated by Cabrerizo et.al for almost contact manifolds [9]. They defined these submanifolds as follows:\\

\noindent {\bf{Definition 2.2 [9].}} {\it{A submanifold $M$ of an almost contact manifold $\bar M$ is said to be a semi-slant submanifold if there exist two orthogonal distributions $D$ and $D_\theta$ such that}}
\begin{enumerate}
\item [(i)] $TM=D\oplus D_\theta\oplus \langle\xi\rangle$
\item [(ii)] {\it{$D$ is an invariant  i.e., $\phi D\subseteq T M $.}}
\item [(iii)] {\it{$D_\theta$ is a slant distribution with slant angle $\theta\neq \frac{\pi}{2}$.}}
\end{enumerate}

\parindent=8mm In the above definition, if $\theta=\pi/2$ then $M$ is contact CR-submanifold of $\bar M$. If $\nu$ is the invariant subspace of the normal bundle $T^\perp M$, then in case of semi-slant submanifolds, the normal bundle $T^\perp M$ can
be decomposed as follows
$$T^\perp M= FD_\theta\oplus \nu.\eqno(2.13)$$

\parindent=8mm
For differential function $\psi$ on $M$, the gradient $grad\psi$ and the Laplacian $\Delta \psi$ of $\psi$ are defined respectively by
$$g(grad \psi, X) = X\psi,\eqno (2.14)$$
$$\Delta \psi =\sum_{i=1}^{n}((\nabla_{e_i}e_i)\psi- e_ie_i \psi), \eqno (2.15)$$
for any vector field $X$ tangent to $M$, where $\nabla$ denotes the Riemannian connection on $M$.

\section{Warped product submanifolds}
\parindent=8mm In this section, we study warped product submanifolds of nearly trans-Sasakian manifolds. We recall the following results on warped products for later use.\\

\noindent
{\bf {Lemma 3.1 [21].}} { Let $M = N_1\times_{f}N_2$ be a warped product manifold with the warping function f. Then
\begin{enumerate}
\item [(i)] $\nabla_XY\in\Gamma(TN_1)$
\item[(ii)] $\nabla_XZ=\nabla_ZX=(X\ln f)Z$
\item[(iii)] $\nabla_ZW ={\nabla_Z}^{N_2}W-(g(Z, W)/f) grad f,$
\end{enumerate}
for any $X, Y\in\Gamma(TN_1)$ and $Z, W\in\Gamma(TN_2)$, where $\nabla$ and $\nabla ^{N_2}$ denote the Levi-Civita connections on $M$ and $N_2$, respectively and $grad f$ is the gradient of $f$.\\

\noindent
{\bf {Corollary 3.1 [21].}} {\it{On a warped product manifold $M = N_1 \times _{f}N_2$, we have}}
\begin{enumerate}
\item [(i)] {\it{$N_1$ is totally geodesic in $M$}}
\item [(ii)] {\it{$N_2$ is totally umbilical in $M.$}}
\end{enumerate}

\parindent=8mm
In the following, we prove the non-existence of warped products of the form $M=N_1\times _fN_2$ in a nearly trans-Sasakian manifold such that $\xi$ is tangent to $N_2$.\\

\noindent
{\bf{Theorem 3.1.}} {\it Let $\bar M$ be a nearly trans-Sasakian manifold that is not nearly Sasakian and $M=N_1\times _fN_2$ be a warped product submanifold of  $\bar M$ such that $\xi$ is tangent to $N_2$, then $M$ is simply a Riemannian product of $N_1$ and $N_2$, where $N_1$ and $N_2$ are any Riemannian submanifolds of $\bar M$.}\\

\noindent
{\it Proof.} For any $X\in \Gamma (TN_1)$, we have $(\bar \nabla_X\phi )\xi+(\bar\nabla_\xi \phi)X=-\alpha X-\beta \phi X.$ By a simple calculation, this relation gives
$$-\phi\bar\nabla_X\xi+\bar\nabla_\xi\phi X-\phi \bar \nabla_\xi X=-\alpha X-\beta\phi X.\eqno(3.1)$$
Taking the inner product with $\phi X$ in (3.1) and using the fact that $\xi$ is tangent to $N_2$, the above equation takes the form
$$g(\bar \nabla_\xi \phi X, \phi X)= - \beta \|X\|^2.\eqno(3.2)$$
Now, since $X$ and $\xi$ are orthogonal vectors, then by (1.2) we can write $\xi g(\phi X, \phi X)= \xi g(X,X),$ which is equivalent to $$g(\bar\nabla_\xi\phi X, \phi X)=g(\bar\nabla_\xi X, X).$$
In view of Lemma 3.1 (ii), the right hand side of the above equation vanishes, hence
$$g(\bar\nabla_\xi\phi X, \phi X)=0.\eqno(3.3)$$
From (3.2) and (3.3), we get$\beta \|X\|^2=0,$ this means that the first factor of the warped product vanishes, which proves the theorem completely.$~\blacksquare$

\parindent=8mm
In view of the above theorem we get a non-existence result about the warped product semi-slant submanifolds in a nearly trans-Sasakian manifold, i.e. there do not exist warped product semi-slant submanifolds $N_\theta \times _fN_T$ and $N_T\times _fN_\theta$ of a nearly trans-Sasakian manifold when the characteristic vector field $\xi $ is a tangent to the second factor. Now we are going to show that the warped product $N_\theta \times _fN_T$ is also a Riemannian product if $\xi$ is tangent to the first factor.\\

\noindent
{\bf{Theorem 3.2.}} {\it There do not exist warped product semi-slant submanifolds of the type $M=N_\theta \times _fN_T$ of a nearly trans-Sasakian manifold $\bar M$ such that $\xi$ is tangent to $N_\theta$, unless $\bar M$ is nearly $\beta$-Kenmotsu.}\\

\noindent
{\it Proof.} Consider $X$ as an arbitrary tangent vector to $N_T$, then making use of (1.3) it follows $(\bar \nabla_X\phi ) \xi + ( \bar \nabla_\xi \phi ) X=  - \alpha X - \beta \phi X.$ This relation can be simplified as
$$- \phi \bar \nabla_X\xi + \bar \nabla_\xi \phi X - \phi \bar \nabla_\xi X =  - \alpha X - \beta \phi X.\eqno(3.4)$$
Taking the inner product with $X$ in (3.4), we get
$$- g(\phi \bar \nabla_X\xi, X) + g(\bar \nabla_\xi \phi X, X) - g(\phi \bar \nabla_\xi X, X) =  - \alpha \|X\|^2.\eqno(3.5)$$
By orthogonality of $X$ and $\xi$ and Lemma 3.1 (ii), the left hand side of (3.5) vanishes identically, hence we reach $\alpha \|X\|^2=0,$ this means that the first factor of the warped product $N_\theta \times _fN_T$ vanishes, which proves the theorem.$~\blacksquare$

\parindent=8mm
From the above discussion, we conclude that there do not exist warped product semi-slant submanifolds of type  $N_\theta \times _fN_T$ in a nearly trans-Sasakian manifold $\bar M$ in both the cases either $\xi$ is tangent to the first factor or to the second. Also, the warped product $N_T\times _fN_\theta$ is just a Riemannian product when the characteristic vector field $\xi$ is tangent to $N_\theta$. Now, we discuss the warped product submanifolds $N_T\times _fN_\theta$ such that $\xi$ is tangent to $N_T$.

\parindent=8mm First, we prove a key lemma characterizing geometric properties of warped product submanifolds $N_T\times _fN_\theta$ of a nearly trans-Sasakian manifold $\bar M$.\\

\noindent
{\bf{Lemma 3.2.}} {\it{Let $M=N_T\times _fN_\theta$ be a warped product semi-slant submanifold of a nearly trans-Sasakian manifold $\bar M$ such that $\xi$ is tangent to $N_T$. Then the following hold}}
\begin{enumerate}
\item [(i)] $\xi \ln f=\beta$,
\item [(ii)] $g(h(X,Y), FZ)=0$,
\item [(iii)] $g(h(\xi, Z), FW)= -\alpha g( Z,W),$
\item [(iv)] $g(h(X,Z), FZ)=-\{(\phi X\ln f)+\alpha \eta (X)\}\|Z\|^2,$
\item [(v)] $g(h(X,Z), FPZ)=-g(h(X,PZ), FZ)=\frac{1}{3}\cos ^2\theta\{(X\ln f)-\beta \eta (X)\} \|Z\|^2,$
\item [(vi)] $g(h(X,X),\zeta)=-g(h(\phi X, \phi X), \zeta)$
\end{enumerate}
{\it{for any $X, Y\in\Gamma(TN_T)$, $Z, W\in\Gamma(TN_\theta)$  and $\zeta\in\Gamma(\nu)$.}}\\

\noindent
{\it Proof.} The first three parts can be proved by the same way as we have proved for contact CR-warped products in [19]. Now, as we consider $\xi$ is tangent to $N_T$, then  for any $X\in\Gamma(TN_T)$ and $Z\in \Gamma(TN_\theta)$, we have
$$(\bar\nabla_X\phi)Z+(\bar\nabla_X\phi)Z=-\alpha\eta(X)Z-\beta\eta(X)\phi Z.$$
Taking the inner product with $Z$, we obtain
$$(\bar\nabla_X\phi)Z+(\bar\nabla_X\phi)Z, Z)=-\alpha\eta(X)\|Z\|^2.\eqno(3.6)$$

Also, we have
$$(\bar\nabla_X\phi)Z=\bar\nabla_X\phi Z-\phi \bar\nabla_XZ.~~~~~~~~~$$
$$~~~~~~~~~~~~~~~~~~~=\nabla_XPZ+h(X, PZ)-A_{FZ}X$$
$$~~~~~~~~~~~~~~~~~~~+\nabla^\perp_XFZ-\phi \nabla_XZ-\phi h(X, Z).$$
Taking the inner product with $Z$ and using Lemma 3.1 (ii), we obtain
$$g((\bar\nabla_X\phi)Z, Z)=0.\eqno(3.7)$$
Similarly, we can obtain
$$g((\bar\nabla_Z\phi)X, Z)=(\phi X\ln f)\|Z\|^2+g(h(X, Z), FZ).\eqno(3.8)$$
Then from (3.6), (3.7) and (3.8), we obtain part (iv) of the lemma. Now, from the structure  of nearly trans-Sasakian manifolds and Lemma 3.1 (ii), we have
$$g((\bar\nabla_X\phi)PZ+(\bar\nabla_{PZ}\phi)X, Z)=\beta\eta(X)\cos^2\theta\|Z\|^2\eqno(3.9)$$
for any $X\in\Gamma(TN_T)$ and $Z\in \Gamma(TN_\theta)$ such that $\xi$ is tangent to $N_T$. Again, by Lemma 3.1 (ii) and Gauss-Weingarten formulas, we obtain
$$g((\bar\nabla_X\phi)PZ, Z)=g(h(X, PZ), FZ)-g(h(X, Z), FPZ)\eqno(3.10)$$
and 
$$g((\bar\nabla_{PZ}\phi)X, Z)=g(h(X, PZ), FZ)+(X\ln f)\cos^2\theta\|Z\|^2.\eqno(3.11)$$
Thus from (3.9), (3.10) and (3.11), we derive
$$2g(h(X, PZ), FZ)-g(h(X, Z), FPZ)=\{\beta\eta(X)-(X\ln f)\}\cos^2\theta\|Z\|^2.\eqno(3.12)$$
Interchanging $Z$ by $PZ$ in (3.12), we obtain
$$-2g(h(X, PZ), FZ)+g(h(X, Z), FPZ)=\{\beta\eta(X)-(X\ln f)\}\cos^2\theta\|Z\|^2.\eqno(3.13)$$
Then, by (3.12) and (3.13), we get
$$g(h(X, PZ), FZ)=-g(h(X, Z), FPZ)\eqno(3.14)$$
which is first equality of the fifth part of the Lemma. The second equality of (v) follows from (3.12) and (3.14). For the last part of the lemma, for any $X\in\Gamma(TN_T)$, we have $\bar \nabla_X\phi X-\phi \bar\nabla_XX=\alpha \|X\|^2 \xi-\eta(X)X-\beta\eta(X)\phi X.$ By means of (2.1), this relation reduces to
$$\nabla_X\phi X+h(\phi X,X)-\phi \nabla_XX-\phi h(X,X)=\alpha\|X\|^2\xi-\eta (X) X- \beta\eta(X)\phi X.$$
Taking the inner product in the above equation with $\phi \zeta$, for any vector $\zeta\in\Gamma(\nu)$, we deduce that
$$g(h(\phi X,X), \phi \zeta)-g(h(X,X), \zeta)=0.\eqno(3.15)$$
Interchanging $X$ by $\phi X$ in the above equation and making use of (1.1) and the fact that $\nu$ is an invariant normal subbundle of $T^\perp M$, it yields
$$-g(h(X, \phi X), \phi \zeta)+ \eta (X) g(h(\xi , \phi X), \phi \zeta)=g(h(\phi X, \phi X), \zeta).\eqno (3.16)$$
Now, by means of (1.3), we derive
$$\nabla_\xi \phi X+ h(\phi X, \xi )- \phi \nabla_\xi X- 2 \phi h(X, \xi )- \phi \nabla_X\xi = \alpha (\eta (X) \xi -X) - \beta \phi X.\eqno(3.17)$$
Taking the inner product with $\phi \zeta$ in (3.17), it follows
$$g(h(\phi X, \xi ) , \phi \zeta )- 2 g(h(X, \xi ), \zeta )=0.$$
Interchanging $\zeta$ by $\phi \zeta$ in the first step and $X$ by $\phi X$ in the second, taking in consideration that $h(\xi, \xi)=0$, we obtain the following couple of tensorial relations
$$g(h(\phi X, \xi ), \zeta ) + 2 g(h(X, \xi ), \phi \zeta )=0,\eqno(3.18)$$
and
$$g(h(X, \xi ), \phi \zeta )+ 2 g(h(\phi X, \xi ) , \zeta )=0. \eqno(3.19)$$
From (3.18) and (3.19), we deduce that
$$g(h(X, \xi ), \phi \zeta )=  g(h(\phi X, \xi ) , \zeta ). \eqno(3.20)$$
In view of (3.19) and (3.20), we get $g(h(X, \xi ), \phi \zeta )= 0.$ Again, interchanging $X$ by $\phi X$ in this relation, it yields
$$g(h(\phi X, \xi ), \phi \zeta )=0.\eqno(3.21)$$
Then  by (3.16) and (3.21), we reach
$$-g(h(\phi X,X), \phi \zeta)-g(h(\phi X,\phi X), \zeta)=0.\eqno(3.22)$$
Thus from (3.15) and (3.22), we get the assertion.$~\blacksquare$

\section{An inequality for warped product submanifolds $N_T\times _fN_\theta$}

In the setting of almost contact structures, many authors have proved  general inequalities in terms of the squared norm of the second fundamental form and the gradient of the warping function in various structures [16, 18, 19, 24]. In fact, all these inequalities are the extension of the original inequality constructed by Chen in the almost Hermitian setting [10]. However, no one proved this relation for warped product semi-slant submanifolds. For this reason, our inequality generalizes the inequalities obtained for CR-warped products in almost contact setting. Another reason is that a nearly trans-Sasakian structure includes all almost contact structures as a special case.

\sloppy
\parindent=8mm
From now on, we shall follow the following orthonormal basis frame of the ambient manifold $\bar M$ for the warped product semi-slant submanifold $M=N_T\times _fN_\theta$ such that $\xi$ is tangent to $N_T$. We shall denote $D$ and $D_\theta$ for the tangent spaces of $N_T$ and $N_\theta$, respectively instead of $TN_T$ and $TN_\theta$. We set $\{e_1,\cdots,e_s, e_{s+1}=\phi e_1,\cdots, e_{(n_1-1=2s)}=\phi e_s, e_{(n_1=2s+1)}=\xi,e_{n_1+1}=e_1^\star,\cdots, e_{n_1+q}=e_q^\star,e_{n_1+q+1}=e_{q+1}^\star=\sec\theta Pe_1^\star,\cdots,e_{(n=n_1+n_2)}=e_{(n_2=2q)}^\star=\sec\theta Pe_q^\star,e_{n+1}=\csc\theta Fe_1^\star,\cdots , e_{n+n_2}=\csc \theta Fe_{n_2}^\star, e_{n+n_2+1}=\bar e_1,\cdots , e_{2m+1}=\bar e_{2l} \}$ as a basis frame of $T\bar M$, then $\{e_1, \cdots , e_s, e_{s+1} =\phi e_1,\cdots, e_{n_1-1}=\phi e_s, e_{n_1}=\xi, e_{n_1+1}=e_1^\star,\cdots, e_{n_1+q}=e_q^\star,e_{n_1+q+1}=e_{q+1}^\star=\sec \theta Pe_1^\star,\cdots,e_{(n=n_1+n_2)}=e_{(n_2=2q)}^\star =\sec\theta Pe_q^\star\}$ are the basis of $TM$,  such that $ e_1, \cdots , e_s, e_{s+1} =\phi e_1,\cdots, e_{n_1-1}=\phi e_s, e_{n_1}=\xi $ are tangent to $D$ and  $e_1^\star,\cdots, e_q^\star, e_{q+1}^\star=\sec \theta Pe_1^\star,\cdots, e_{(n_2=2q)}^\star =\sec \theta Pe_q^\star$ are tangent to $D_\theta$, hence $\{ e_{n+1}=\csc \theta Fe_1^\star,\cdots, e_{n+n_2}=\csc \theta Fe_{n_2}^\star, e_{n+n_2+1}=\bar e_1, \cdots, e_{2m+1}=\bar e_{2l}\}$ are the basis of the normal bundle $T^\perp M$, such that $e_{n+1}=\csc\theta Fe_1^\star,\cdots, e_{n+n_2}=\csc\theta Fe_{n_2}^\star$ are tangent to $FD_\theta$ and $ e_{n+n_2+1}=\bar e_1,\cdots, e_{2m+1}=\bar e_{2l}$ are tangent to the invariant normal subbundle $\nu$ with dimension $2l$. We use this frame in the following theorem.\\

\noindent
{\bf{Theorem 4.1.}} {\it{Let  $M=N_T\times{_{f}N_\theta}$ be a warped product semi-slant submanifold of a nearly trans-Sasakian  manifold $\bar M$ such that $\xi$ is tangent to  $N_T$, where $N_T$ and $N_\theta$ are invariant and proper slant submanifolds of dimensions $2s+1$ and $2q$, respectively. Then,}}
\begin{enumerate}
\item [(i)] {\it The second fundamental form of $M$ satisfies the following inequality}
$$\|h\|^2\geq 2q[\{\frac{2}{9}\cot^2\theta + 2\csc^2\theta\}(\|grad(\ln f)\|^2-\beta^2)+\alpha^2].\eqno(4.1)$$
\item [(ii)] {\it{If the equality sign in (i) holds identically, then $N_T$ and $N_\theta $ are totally geodesic and totally umbilical submanifolds in $\bar M$, respectively.}}
\end{enumerate}

\noindent
{\it {Proof.}} In view of the above frame and the definition of the second fundamental form, it is straightforward to get the following expansion
$$\|h\|^2=\sum_{r=n+1}^{2m+1} \sum_{i, j=1}^{n} g(h(e_i, e_j), e_r)^2~~~~~~~~~~~~~~~~~~~~~~~~~~~$$
$$~~~~~~~~~~~~~~~~~~~~~~=\sum_{r=n+1}^{n+n_2} \sum_{i, j=1}^{n} g(h(e_i, e_j), e_r)^2+\sum_{r=n+n_2+1}^{2m+1} \sum_{i, j=1}^{n} g(h(e_i, e_j), e_r)^2$$
$$~~~~~~~~~~~~~~~~~~~~~~\geq\sum_{r=n+1}^{n+n_2} \sum_{i, j=1}^{n} g(h(e_i, e_j), e_r)^2=\sum_{l=2s+2}^{n} \sum_{i, j=1}^{n} g(h(e_i, e_j), \phi e_l)^2.$$
Using the orthonormal frame of $D$ and $D_\theta$, it follows
$$\|h\|^2\geq\sum_{l=2s+2}^{n} \sum_{i, j=1}^{2s+1} g(h(e_i, e_j), \phi e_l)^2+2\sum_{j,l=2s+2}^{n}\sum_{i=1}^{2s+1}  g(h(e_i, e_j), \phi e_l)^2$$
$$+\sum_{i, j, l=2s+2}^{n}  g(h(e_i, e_j), \phi e_l)^2.~~~~~~~~~~~~\eqno(4.2)$$
By Lemma 3.2 (ii), the first term of the right hand side in (4.2) is identically zero,
so let us compute the next term
$$\|h\|^2\ge 2\sum_{j, l=2s+2}^{n}\sum_{i=1}^{2s} g(h(e_i, e_j), \phi e_l)^2+2\sum_{j,l=2s+2}^{n} g(h(\xi, e_j), \phi e_l)^2.\eqno(4.3)$$
Making use of Lemma 3.2 (iii), the second term of the right hand side in (4.3) can be evaluated, while by means of the orthonormal frame the first term is expanded to give four terms, as a result (4.3) takes the following form
$$\|h\|^2\ge 2\csc^2 \theta\sum_{j=1}^{q}\sum_{i=1}^{2s}g(h(e_i,e_j), Fe_j^\star)^2~~~~~~~~~~~~~~~~$$
$$+2\csc^2\theta \sec^2\theta\sum_{j=1}^{q}\sum_{i=1}^{2s}  g(h(e_i,Pe_j^\star), Fe_j^\star)^2$$
$$+2\csc^2\theta \sec^2\theta\sum_{j=1}^{q}\sum_{i=1}^{2s}  g(h(e_i,e_j), FPe_j^\star)^2$$
$$~~~+2\csc^2\theta \sec^4\theta\sum_{j=1}^{q}\sum_{i=1}^{2s}  g(h(e_i,Pe_j^\star), FPe_j^\star)^2$$
$$+2\sum_{j, l=2s+2}^{n} (-\alpha g(e_j,e_l))^2.~~~~~~~~~~~~~~~~~~~~\eqno(4.4)$$
Lemma 3.2 (iii)-(v) can be used to substitute the inner products in the right hand side of (4.4), thus it yields
$$\|h\|^2\ge 2\csc^2\theta\sum_{j=1}^{q}\sum_{i=1}^{2s}((\phi e_i\ln f)+\alpha \eta(e_i))^2 \|e_j\|^4~~~~~~~~~~~~~$$
$$~+\frac{2}{9}\cos ^2\theta \csc^2\theta\sum_{j=1}^{q}\sum_{i=1}^{2s}((e_i\ln f)-\beta \eta(e_i))^2\|e_j\|^4$$
$$~+\frac{2}{9}\cos ^2\theta \csc^2\theta\sum_{j=1}^{q}\sum_{i=1}^{2s}((e_i\ln f)-\beta \eta (e_i))^2\|e_j\|^4$$
$$~~~~+2\csc^2\theta\sum_{j=1}^{q}\sum_{i=1}^{2s}((\phi e_i\ln f)+\alpha \eta (e_i)^2\|e_j\|^4+2q\alpha^2.\eqno(4.5) $$
In view of the assumed orthonormal frame, the $1-$form $\eta (e_i)$ is identically zero for all $i\in \{1,\cdots, 2s\}$, hence we reach
$$\|h\|^2\geq 4\csc^2 \theta\sum_{j=1}^{q}\sum_{i=1}^{2s}  (\phi e_i\ln f)^2 \|e_j\|^4~~~~~~~~~~~~~~~$$
$$~~~~~~~~~~~~+\frac{4}{9}\csc^2\theta \cos^2\theta\sum_{j=1}^{q}\sum_{i=1}^{2s}(e_i\ln f)^2 \|e_j\|^4+2q\alpha^2 .\eqno(4.6)$$
Then from (2.14) and Lemma 3.2 (i), the above inequality takes the form
$$\|h\|^2\ge 2q[\{\frac{2}{9}\cot^2\theta + 2\csc^2\theta\}(\|\nabla \ln f\|^2-\beta^2)+\alpha^2],$$
which is the inequality (i). Now, assume that the equality sign in (4.1) holds identically, then from (4.2), (4.3) and Lemma 3.2 (ii), we deduce that
$$h(D, D)=0,~~~h(D_\theta, D_\theta)=0,~~~h(D, D_\theta)\subset FD_\theta .\eqno(4.7)$$
Hence, combine statement of Corollary 3.1 (i) with the first condition in (4.7) shows that $N_T$ is totally geodesic in $\bar M$.
On the other hand, if we denote by $h^\theta$ the second fundamental form of  $N_\theta$ in
$M$. Then, we get
$$g(h^\theta (Z, W) ,X)=g(\nabla_ZW, X)=-(X\ln f)g(Z,W)=-g(Z,W)g(\nabla \ln f, X),$$
which is equivalent to
$$h^\theta (Z,W)= -\nabla\ln f g(Z,W).\eqno(4.8)$$
This means that $N_\theta$ is totally umbilical in $M$, thus the second condition of (4.7) with (4.8) and Corollary 3.1 (ii) imply that $N_\theta$ is totally umbilical in $\bar M$. Also, all three conditions of (4.7) gives the minimality of $M.~\blacksquare$\\

\parindent=0mm
{\bf{Note.}} In the inequality (5.1), if $\alpha=0$ and $\beta=1$, then it reduces to
$$\|h\|^2\geq 2q[\{\frac{2}{9}\cot^2\theta + 2\csc^2\theta\}(\|\nabla\ln f\|^2-1)]$$
which is the inequality for nearly Kenmotsu manifolds. Also, If $\alpha=1$ and $\beta=0$, then the inequality reduces for the nearly Sasakian manifolds. The equality cases can also be discussed.\\

\noindent
{\bf {Remark 1.}} {\it Theorem 3.1 in  [16], Theorem 3.4 in [18] and Theorem 3.2 in [24] are the special cases of the above inequality.}\\

\noindent
{\bf {Remark 2.}} {\it{The above inequality generalizes Theorem 4.1 in [19]}}.

\section{Another inequality for warped products}

\noindent
Let $\varphi :M=N_1\times _fN_2 \longrightarrow \bar M$ be an isometric immersion of a warped product $N_1\times _fN_2$ into a Riemannian manifold $\bar {M}$ of constant sectional curvature $c$. Denote by $n_1,~n_2,~n$ the dimensions of $N_1,~N_2,~N_1\times_f N_2$, respectively. Then for unit vector fields $X,~Z$ tangent to $N_1,	N_2$, respectively, we have
$$K(X\wedge Z)=g(\nabla_Z\nabla_XX-\nabla_X\nabla_ZX, Z)$$
$$=(1/f)\{(\nabla_XX)f-X^2f\}.\eqno(5.1)$$
If we choose the local orthonormal frame $e_1,\cdots,e_n$ such that $e_1,\cdots, e_{n_1}$ are tangent to $N_1$ and $e_{n_1+1},\cdots,e_n$ are tangent to $N_2$, then we have
$$\frac{\Delta f}{f}=\sum_{i=1}^{n_1}K(e_i\wedge e_j)\eqno(5.2)$$
for each $j=n_1+1,\cdots, n$.

\parindent=8mm
In this section, our aim is to develop a new method which is giving a useful formula for the squared norm of the mean curvature vector $\vec H$ under $\varphi$, Geometrically, this formula declares the $N_T$-minimality of $\varphi$. 

\parindent=8mm
We know that 
$$\|H\|^2=\frac {1}{n^2} \sum_{r=n+1}^{2m+1}(h_{11}^r+ \cdots + h_{nn}^r)^2.$$
Taking in consideration that ($n=n_1+n_2$), where $n_1$ and $n_2$ are the dimensions of $N_T$ and $N_\theta$, respectively, then it follows
$$\|H\|^2=\frac {1}{n^2} \sum_{r=n+1}^{2m+1}(h_{11}^r+ \cdots +h_{n_1n_1}^r+ h_{n_1+1 n_1+1}^r+ \cdots + h_{nn}^r)^2.$$
Moreover, for every $r\in \{n+1, \cdots , 2m+1\}$, using the frame of $D$ and the fact that $h(\xi, \xi )=0$, then $n_1$ coefficients of the right hand side can be decomposed as
$$(h_{11}^r+\cdots +h_{n_1n_1}^r+ h_{n_1+1 n_1+1}^r+ \cdots + h_{nn}^r)^2~~~~~~~~~~~~~~~~~~~~~~~~~~~~~~~~~~~~~~~~~~$$
$$~~~~~~~~=(h_{11}^r+ \cdots +h_{ss}^r+h_{s+1s+1}^r+ \cdots + h_{2s2s}^r+h_{\xi \xi}^r+ h_{n_1+1 n_1+1}^r+ \cdots + h_{nn}^r)^2$$
$$~~~~~~~=(h_{11}^r+ \cdots +h_{ss}^r+h_{s+1s+1}^r+ \cdots + h_{2s2s}^r+ h_{n_1+1 n_1+1}^r+ \cdots + h_{nn}^r)^2.\eqno (5.3)$$
From (2.7), we know that $e_r$ belongs to the normal bundle $TM^\perp$, for every $r\in \{n+1, \cdots , 2m+1\}$, then in view of (2.13) we have two cases, either it belongs to $FD_\theta$ or to $\nu$.

\parindent=8mm
{\it{Case (i).}} If $e_r\in\Gamma(FD_\theta)$, then from Lemma 3.2 (ii), we know that $g(h(X,X), FZ)=0$, for any $X\in\Gamma(D)$ and $Z\in\Gamma(D_\theta)$, consequently (5.3) reduces to
$$(h_{11}^r+ \cdots +h_{n_1n_1}^r+ h_{n_1+1 n_1+1}^r+ \cdots + h_{nn}^r)^2=( h_{n_1+1 n_1+1}^r+ \cdots + h_{nn}^r)^2.\eqno(5.4)$$

\parindent=8mm
{\it{Case (ii).}} If $e_r\in \Gamma(\nu)$, then by means of Lemma 3.2 (vi), we can make an expansion of (5.3) as 
$$(h_{11}^r+ \cdots +h_{n_1n_1}^r+ h_{n_1+1 n_1+1}^r+ \cdots + h_{nn}^r)^2~~~~~~~~~~~~~~~~~~~~~~~~~~~~~~~~~~~~~~~~~~~~~~~~~~~~~~~~~~~~~~~~~~~~~~$$
$$~~=(g(h(e_1,e_1),e_r)+ \cdots +g(h(e_s, e_s), e_r)+g(h(\phi e_1, \phi e_1), e_r)$$
$$+\cdots+g(h(\phi e_s, \phi e_s), e_r)+h_{n_1+1 n_1+1}^r+ \cdots + h_{nn}^r)^2~~~~~~$$
$$~~~~~~=(g(h(e_1,e_1),e_r)+ \cdots +g(h(e_s, e_s), e_r)-g(h(e_1,e_1),e_r)- \cdots$$
$$~~~~-g(h(e_s, e_s), e_r)+ h_{n_1+1 n_1+1}^r+ \cdots + h_{nn}^r)^2~~~~~~~~~~~~~~~~~~~~$$
$$=( h_{n_1+1 n_1+1}^r+ \cdots+h_{nn}^r)^2.~~~~~~~~~~~~~~~~~~~~~~~~~~~~~~~~~~~~~~~\eqno(5.5)$$
Then from (5.4) and (5.5), we can deduce that
$$(h_{11}^r+ \cdots +h_{n_1n_1}^r+ h_{n_1+1 n_1+1}^r+ \cdots + h_{nn}^r)^2=( h_{n_1+1 n_1+1}^r+ \cdots + h_{nn}^r)^2$$
for every normal vector $e_r$ belongs to the normal bundle $T^\perp M$. In other words
$$\sum_{r=n+1}^{2m+1}(h_{11}^r+ \cdots +h_{n_1n_1}^r+ h_{n_1+1 n_1+1}^r+\cdots + h_{nn}^r)^2=\sum_{r=n+1}^{2m+1}( h_{n_1+1 n_1+1}^r+ \cdots + h_{nn}^r)^2.$$
By the end of this discussion, we can state the following lemma.\\

\noindent
{\bf {Lemma 5.1.}} { Let $\varphi :M=N_T\times _fN_\theta \longrightarrow \bar M$ be an isometric immersion from a warped product semi-slant submanifold into a nearly trans-Sasakian manifold $\bar {M}$. Then, we have

$$\|\vec H\|^2=\frac{1}{n^2}\sum_{r=n+1}^{2m+1}( h_{n_1+1 n_1+1}^r+\cdots+ h_{nn}^r)^2,$$
i.e., $\varphi$ is $N_T$-minimal immersion, where $\vec H$ is the mean curvature vector and $n_1$, $n_2$, $n$ and ($2m+1$) are the dimensions of $N_T$, $N_\theta$, $M$ and $\bar M$, respectively.

\parindent=8mm
From the Gauss equation and the above key Lemma 5.1, we are able to state and prove the following general inequality.\\

\noindent
{\bf {Theorem 5.1.}} {\it{Let $\varphi :M=N_T\times _fN_\theta \longrightarrow \bar M$ be an isometric immersion from a  warped product semi-slant submanifold into a nearly trans-Sasakian manifold $\bar {M}$ such that $\xi$ is tangent to $N_T$. Then, we have}}
\begin{enumerate}
\item [(i)] $\frac{1}{2}\|h\|^2\ge \bar \tau (TM)-\bar \tau (TN_T)-\bar \tau (TN_\theta)-\frac {n_2 \Delta f}{f},$

\noindent {\it{where $n_2$ is the dimension of $N_\theta$.}}
\item [(ii)] {\it{If the equality sign in (i) holds identically, then $N_T$ and $N_\theta$ are totally geodesic and totally umbilical submanifolds in $\bar M$, respectively.}}
\end{enumerate}

\noindent
{\it{Proof.}} { We start by recalling (2.9) as a consequence of (2.5) as 
$$\|h\|^2= -2 \tau+2 \bar\tau (TM)+ n^2\|H\|^2.$$
Making use of (2.6) in the above equation, we deduce
$$\|h\|^2= -2 \sum_{i=1}^{n_1} \sum_{j=n_1+1}^{n} K(e_i\wedge e_j) - 2 \tau (T N_T)-2 \tau (T N_\theta )+2 \bar\tau(TM)+ n^2\|H\|^2.$$
Then by Lemma 3.2 and the relation (2.8), it follows
$$\|h\|^2=- \frac{ 2n_2\Delta f}{f}-2\bar \tau(T N_T)-2\sum_{r=n+1}^{2m+1}\sum_{1\le i<k\le n_1} (h_{ii}^r h_{kk}^r- (h_{ik}^r)^2)-2\bar\tau (TN_\theta )$$
$$-2\sum_{r=n+1}^{2m+1}\sum_{n_1+1\le j<t\le n} (h_{jj}^r h_{tt}^r- (h_{jt}^r)^2)+2 \bar \tau (TM)+ n^2\|H\|^2.$$
The above equation is equivalent to the following form
$$\|h\|^2=- \frac{ 2n_2\Delta f}{f} - 2\bar\tau (TN_T)-\sum_{r=n+1}^{2m+1}\sum_{1\le i\ne k\le n_1} (h_{ii}^r h_{kk}^r- (h_{ik}^r)^2)-2 \bar\tau (TN_\theta )$$
$$+2 \bar \tau (TM)-\sum_{r=n+1}^{2m+1}\sum_{n_1+1\le j\ne t\le n} (h_{jj}^r h_{tt}^r- (h_{jt}^r)^2)+ n^2\|H\|^2.$$
The above equation takes the following form when we add and subtract the same term on the right hand side
$$\|h\|^2=-\frac{ 2n_2\Delta f}{f}-2\bar \tau (TN_T)-\sum_{r=n+1}^{2m+1} ((h_{11}^r)^2+\cdots+ (h_{n_1n_1}^r)^2)~~~~~~~~~~~~~~~~~~~~~~~~~~$$
$$-\sum_{r=n+1}^{2m+1}\sum_{1\le i\ne k\le n_1} (h_{ii}^r h_{kk}^r- (h_{ik}^r)^2)+\sum_{r=n+1}^{2m+1} ((h_{11}^r)^2 +\cdots + (h_{n_1n_1}^r)^2)~~$$
$$-2 \bar\tau(TN_\theta )-\sum_{r=n+1}^{2m+1}\sum_{n_1+1\le j\ne t\le n} (h_{jj}^r h_{tt}^r- (h_{jt}^r)^2)+2 \bar \tau (TM)+ n^2\|H\|^2,$$
$$~~~~=-\frac{ 2n_2\Delta f}{f}-2\bar \tau (TN_T)+\sum_{r=n+1}^{2m+1}\sum_{i,k=1}^{n_1}(h_{ik}^r)^2- \sum_{r=n+1}^{2m+1} (h_{11}^r+\cdots+h_{n_1n_1}^r)^2$$
$$~-2\bar\tau (TN_\theta )-\sum_{r=n+1}^{2m+1}\sum_{n_1+1\le j\ne t\le n} (h_{jj}^r h_{tt}^r- (h_{jt}^r)^2)+2 \bar \tau (TM)+ n^2\|H\|^2.$$
Similarly, we can add and subtract the same term for the sixth term in the above equation and finally, we derive
$$\|h\|^2=-\frac{ 2n_2\Delta f}{f}+2 \bar \tau (TM)-2\bar\tau (TN_T)+\sum_{r=n+1}^{2m+1}\sum_{i,k=1}^{n_1}(h_{ik}^r)^2~~~~~~~~~~~~~~~~~~~~~~$$
$$+\sum_{r=n+1}^{2m+1}\sum_{j,t=n_1+1}^{n}(h_{jt}^r)^2-\sum_{r=n+1}^{2m+1} (h_{11}^r+ \cdots+h_{n_1n_1}^r)^2-2 \bar\tau(TN_\theta )$$
$$-\sum_{r=n+1}^{2m+1} (h_{n_1+1n_1+1}^r+\cdots+h_{nn}^r)^2+ n^2\|H\|^2.~~~~~~~~~~~~~~~~~~~~~~~$$
Taking account of Lemma 5.1, we get the inequality (i). For the equality case, from the last relation we get
$$\sum_{r=n+1}^{2m+1}\sum_{i,k=1}^{n_1}g(h(e_i, e_k), e_r)=0\eqno(5.6)$$
and
$$\sum_{r=n+1}^{2m+1}\sum_{j,t=n_1+1}^{n}g(h(e_j, e_t), e_r)=0.\eqno(5.7)$$
From (5.6) and (5.7), we obtain that the immersion $\varphi:M\to \bar M$ is totally geodesic. Also, from Corollary 3.1 we know that the immersion $N_T\to M$ is totally geodesic  and the immersion $N_\theta\to M$ is totally umbilical, hence the result (ii).$~\blacksquare$

\noindent Author's addresses:\\

Abdulqader Mustafa

\noindent Institute of Mathematical Sciences, Faculty of Science,University of Malaya, 50603 Kuala Lumpur, MALAYSIA

\noindent {\it E-mail}: {\tt abdulqader.mustafa@yahoo.com}

\smallskip

Siraj Uddin

\noindent Institute of Mathematical Sciences, Faculty of Science,University of Malaya, 50603 Kuala Lumpur, MALAYSIA

\noindent {\it E-mail}: {\tt siraj.ch@gmail.com}

\smallskip

Bernardine R. Wong

\noindent Institute of Mathematical Sciences, Faculty of Science,University of Malaya, 50603 Kuala Lumpur, MALAYSIA

\noindent {\it E-mail}: {\tt bernardr@um.edu.my}

\end{document}